\documentclass[11pt]{article}
\usepackage{latexsym}
 \usepackage{amssymb}
 \usepackage{amsmath}
 \usepackage{amsthm}
 \usepackage{graphicx}

\newcommand{ \un }{1\!\!1}
 \newcommand{ \p }{\mathbb{P} }
 \newcommand{ \pa }{\mathbb{P}^{\alpha} }

 \newcommand{ \E }{\mathbb{E}}

 \newcommand{ \R }{ \mathbb{R} }
 \newcommand{ \Z }{ \mathbb{Z} }
 \newcommand{\N}{ \mathbb{N} }

 \newcommand{ \f }{ \mathcal{F} }

\newcommand{ \g}{ \mathcal{G} }

\newcommand{ \ct }{ \textrm{const} }

\newcommand{ \G}{ \mathcal{G} }

 \newcommand{ \lo }{ \mathcal{L} }

\newtheorem{The}{{\bf Theorem}}[section]
\theoremstyle{definition}

\theoremstyle{plain}
 
 \newtheorem{Cor}[The]{\bf Corollary}
 
 \theoremstyle{definition}
\newtheorem{Rem}[The]{{\bf Remark}}

\makeatletter

\@addtoreset{equation}{section} \makeatother

\title{Constants of concentration for a simple recurrent random walk on random environment
\\ \vspace{1cm}
 \large{Pierre Andreoletti} \footnote{Laboratoire MAPMO - C.N.R.S. UMR 6628 - F\'ed\'eration Denis-Poisson, Universit\'e d'Orl\'eans, 
(Orl\'eans France). \newline \vspace{0.1cm}  $\quad$  MSC 2000 60G50; 60J55. \newline \vspace{0.5cm} \textit{Key words and phrases :  random environment, random walk, recurrent regime, local time, concentration} } \textrm{ }    }

\begin{document}

\maketitle


\noindent \\ \textbf{Abstract:} We clarify the asymptotic of the limsup of the size of the neighborhood of concentration of Sinai's walk improving the result in \cite{Pierre3}. Also we get the almost sure limit of the number of points visited more than a small but fixed proportion of a given amount of time. 


\section{Introduction, and Results}

The one-dimensional recurrent random walk on random environment (thereafter abbreviated RWRE) we treat here, also called Sinai's walk, has the property to be localized at an instant $n$ in a neighborhood of a point of the lattice. This point can be described from the random environment and depends on $n$ (see \cite{Sinai}, \cite{Golosov1}, \cite{Zeitouni}, \cite{Pierre1}), its limit distribution is also known (see \cite{Golosov}, \cite{Kesten2}). \\
In addition to this aspect of localization, this random walk has the property to spend a large amount of time in a single point of the lattice, this property was first en-lighted by \cite{Revesz}. 
To be more precise, let us define the local time ($\lo$) of the random walk $(X_n,n \in \N)$. Let $k\in \Z$ and $n \in \N^*$
\begin{eqnarray*}
& & \lo(k,n)=\sum_{i=1}^n\un_{X_i=k}, \\ & &  \lo(A,\Z)=\sum_{k \in A}\lo(k,n), \ A\subset \Z, \\
& &  \lo^*(n)=\sup_{k \in \Z}\lo(k,n).
 \end{eqnarray*}
\cite{Revesz} shows that $\limsup_{n}\frac{\lo^*(n)}{n}\geq \ct>0, \p\ a.s.$ and more recently \cite{GanPerShi} get the more precise result $\limsup_{n}\frac{\lo^*(n)}{n} = c_1, \p\ a.s.$ giving an explicit formula for $c_1$ as a function of the support of the distribution of the random environment. \\
In this paper we mainly answer one of the question asked in \cite{Pierre3}, more precisely we are interested in random variables  related to the local time: let $x\in \Z$, $0\leq \beta \leq1$, $r \in \N$, $I_{r}(x)=\{x-r,x-r+1,\cdots,x,\cdots,x+r\} \label{defir}$, we define $$Y_{n,\beta}=\inf_{x\in \Z}\inf(r>0,\lo(I_{r}(x),n)\geq \beta n).$$
$Y_{n,Ê\beta} $ measure the size of a neighborhood where the walk spends more than a proportion $\beta$ of the total amount of time $n$. We call $Y_{n,\beta}$ the concentration variable of $(X_n,n\in \N)$, the concentration is the property of the random walk to spend a large amount of time in a negligible interval comparing to a typical fluctuation of the walk.
The first asymptotics for $Y_{n,Ê\beta}$ can be found in \cite{Pierre2}, more recently \cite{Pierre3} get the following 
\begin{eqnarray}
\liminf_nY_{n,Ê\beta} \leq const(1-\beta)^{-2}, \ \p\ a.s. \label{eq1}
\end{eqnarray}
Some more work was needed to improve the inequality \ref{eq1} into an equality, and especially to show that $(1-\beta)^{-2}$ is not the good behavior. We get it here mainly by applying a method of \cite{GanPerShi}. \\
In this paper we are also interested in a second variable that we define below, let $\delta>0$:
\begin{eqnarray*}
Z_{n,Ê\delta}=\sum_{x \in \Z} \un_{\{\lo(k,n) \geq \delta n\}}.
\end{eqnarray*}
$Z_{n,Ê\delta}$ counts the number of points of the lattice visited more than $\delta n$ times. Our second result give the function $m$ that satisfies for $\delta$ small enough  $\limsup_n Z_{n,Ê\delta}=m(\delta), \p.a.s.$. \\
It is time now to define the model and present the results. The simple random walk on random environment we are dealing with is defined from two processes, the first one, called random environment, is a sequence of i.i.d. random variables $\alpha \equiv (\alpha_i,i \in \Z)$ with distribution $P$, each of the $\alpha_i$ belonging to the interval $(0,1)$. The second one $(X_n, n \in \N)$ is a birth and death process for all fixed $\alpha$,  with transition of probabilities given by $\pa[X_{n+1}=i+1| X_n=i]=1-\pa[X_{n+1}=i-1| X_n=i]=\alpha_i$. 
The whole process is defined on the probabilty space $(\Omega\times\Z^{\N},\f\times \g,\p)$ where $\Omega:=(0,1)^{\Z}$ is the state space of the random enviromnent equiped with its Borel $\sigma$-field $\f$, $\g$ is the Borel $\sigma$-field associated to $\Z^{\N}$ and for all $F \in \f$ and $G \in \G$, $\p[F \times G]=P\otimes \pa[F \times G]=\int_F P(dw_1) \int_G \p^{\alpha(w_1)}(dw_2)$.
The one-dimensional RWRE we are interested in is almost surely recurrent for almost all environment. \cite{Solomon} proved that the necessary and sufficient hypothesis to get such a process is
\begin{eqnarray}
E[\epsilon_0]=0, \label{hyp0}
\end{eqnarray}
 where $(\epsilon_i \equiv \log \frac{1-\alpha_i}{\alpha_i}, i \in  \Z)$. We also add the hypothesis 
 \begin{eqnarray}
 Var[\epsilon_0]>0 \label{hyp1}
 \end{eqnarray}
  to get a non-tivial RWRE, and 
\begin{eqnarray}  
P(| \epsilon_0 |\leq \eta )=1, \eta>0 \label{hyp2}
\end{eqnarray}
for simplicity.  
We can now give the results which depend on the following constants, 
\begin{eqnarray*}
& & \bar{A}:=\sup\{x : x \in supp(\alpha_0)\} \in (1/2,1], \textrm{ and } \bar{\alpha}:=\inf\{x : x \in supp(\alpha_0) \}\in [0,1/2), \\
& & \tilde{\alpha}=\bar{\alpha}/(1-\bar{\alpha}) \textrm{ and } \tilde{A}=(1-\bar{A})/\bar{A}.
\end{eqnarray*}
First le us recall one of the results of \cite{GanPerShi}, and give the explicit formula for $c_1$ that will be useful in the sequel:
\begin{The} (\cite{GanPerShi}) \label{GanPerShi} \label{thgs}Assume that \ref{hyp0}, \ref{hyp1} and \ref{hyp2} are verified, then
$$\ \p \ a.s. \ \liminf_n \frac{\lo^*(n)}{n}=c_1,
 $$
where $c_1=(2\bar{A}-1)(1-2\bar{\alpha})/(2(\bar{A}-\bar{\alpha})\min(\bar{A},1-\bar{\alpha})).$
\end{The} 
\noindent Our first result is the following, 
\begin{The} \label{th1} Assume that \ref{hyp0}, \ref{hyp1} and \ref{hyp2} are verified, for all $0\leq \beta \leq 1$,
\begin{eqnarray}\ \p \ a.s. \ \liminf_n Y_{n,\beta}=  \left\{ \begin{array}{ll} 0 & \textrm{if }\ 0\leq \beta \leq c_1  , \\
f(\beta) & \textrm{if }  c_1 < \beta <1, \\
+\infty &   \textrm{if }  \beta=1,     \end{array} \right. 
\end{eqnarray} 
where $f(\beta)$ is the smallest integer that satisfies
\begin{eqnarray}
& & \frac{(1-\tilde{\alpha})(1-\tilde{A})}{2(1-\tilde{\alpha}\tilde{A})}\sup_{x\in \Z} \left(\sum_{k \in I_{f(\beta)}(x)} F(k)+F(k-1)\right) \geq \beta, \label{eq1.7b} \textrm{with} \\ 
& & F(l)=\left\{ \begin{array}{ll} \tilde{\alpha}^l(\tilde{\alpha}\un_{\{\tilde{\alpha}<\tilde{A}\}}+\un_{\{\tilde{\alpha} \geq\tilde{A}\}}) & \textrm{if}\ l>0, \\
\tilde{\alpha}\un_{\{\tilde{\alpha}<\tilde{A}\}}+\un_{\{\tilde{\alpha} \geq\tilde{A}\}} & \textrm{if }  l=0, \\
   \tilde{A}^{-l-1} (\un_{\{\tilde{\alpha}<\tilde{A}\}}+\tilde{A}\un_{\{\tilde{\alpha} \geq\tilde{A}\}}) & \textrm{if }  l <0,     \end{array} \right.
 \end{eqnarray}
and $\un$ the indicator function. 
\\ Moreover if $\bar{A}=1-\bar{\alpha}$, the result get simpler and $f(\beta)$ is the smallest integer that satisfies: 
\begin{eqnarray}
1-\frac{(1+\tilde{\alpha})\tilde{\alpha}^{f(\beta)-1}}{2} \geq \beta.
\end{eqnarray}
\end{The} 

Notice that $\max(F(k), k \in \Z)=1$, we can achieve this maximum by taking $l=0$ if $\tilde{\alpha} \geq\tilde{A}$ or $l=-1$ if $\tilde{\alpha} <\tilde{A}$. Notice also that the site, say $a^+$, where the supremum in \ref{eq1.7b} is achieved depends on $\tilde{\alpha}$ and $\tilde{A}$: if $\tilde{\alpha} \geq \tilde{A}$ then $a^+ \geq 0$ otherwise $a^+ \leq 0$. For the random walk this means that the center of the interval of concentration in the infinite valley, defined in the next section, is not necessarily the site $0$, and is likely to be in a region where the environment is the flattest. 

The other interesting part of Theorem \ref{th1} is that it leads directly to the asymptotic case when $\beta$ gets close to one:

\begin{Cor} Assume that \ref{hyp0}, \ref{hyp1} and \ref{hyp2} are verified, 
\begin{eqnarray}\ \p \ a.s.\ \lim_{\beta \rightarrow 1^-}\ \liminf_n \frac{Y_{n,\beta}}{|\log(1-\beta)|}=\max\left(\frac{1}{|\log\tilde{A}|}, \frac{1}{|\log\tilde{\alpha}|}\right). \label{eq24} 
\end{eqnarray} 
\end{Cor}

\begin{Rem} 
The upper bound obtained in \cite{Pierre3} (see \ref{eq1}) was pretty far from reality. 
 The case $\beta=1$ can be easily deduce from \ref{eq24}, the case $0 \leq \beta \leq c_1$ is a direct consequence of Theorem \ref{GanPerShi}.
\end{Rem}
\noindent The second theorem deals with $Z_{n,Ê\delta}$
\begin{The} \label{th2} Assume that \ref{hyp0}, \ref{hyp1} and \ref{hyp2} are verified, 
\begin{eqnarray}
 \p\ a.s.\  \lim_{\delta \rightarrow 0^+} \delta \limsup_n Z_{n,Ê\delta} = 1.
\end{eqnarray}
\end{The}
\noindent This therorem shows that when $\delta$ is small, $\limsup_{n}Z_{n,\delta}$ behaves like $1/\delta$, that is we get the same limit behavior than for a simple symmetric random walk with finite state space given by $\{1,2,\cdots,1/\delta\}$ with reflecting barrier at 1 and $1/\delta$.
\begin{Rem} 
Notice that by definition of the random walk $\limsup_n Z_{n,0} =+\infty,\ \p.a.s$. Thanks to Theorem \ref{thgs} for all $\delta>c_1$ $\limsup_n Z_{n,\delta} =0, \p.a.s$, finally regarding to (1.12) in \cite{GanPerShi}, $\liminf_n Z_{n,Ê\delta} =0, \p.a.s$.  
\end{Rem}

\section{Proof of the results}

A fundamental notion for the study of Sinai's random walk is the random potential $(S_k,k\in \Z)$ associated with the random environment:
\begin{eqnarray} 
 && S_{k}:=\left\{ \begin{array}{ll} \sum_{1\leq i \leq k} \epsilon_i, & \textrm{ if }\ k =1,2, \cdots \\
  \sum_{k+1 \leq i \leq 0} - \epsilon_i , &   \textrm{ if }\  k=-1,-2,\cdots \end{array} \right.  \nonumber \\
&& S_{0}:=0, \nonumber
 \end{eqnarray}
recall that $\epsilon_i=\frac{1-\alpha_i}{\alpha_i}$ for all $i \in \Z$. Let $(W^+(x),x \in \Z)$ be a collection of random variables distributed as $S$ conditioned to stay non-negative for all $x>0$ and strictly positive for $x<0$, and in the same way $(W^-(x),x \in \Z)$ is a collection of random variables distributed as $S$ but conditioned to stay strictly positive for all $x>0$ and non-negative for $x<0$. Notice that under \ref{hyp0}, $W^+$ and $W^-$ are well defined and $\sum_{x \in \Z}\exp(-W^{\pm}(x))<\infty, \ a.s.$ (see \cite{Bertoin} or \cite{Golosov}). From $W^+$ and $W^-$ we define two random probability measures, 
$\mu^+$ and $\mu^-$ as follows, let $x \in \Z$, 
\begin{eqnarray}
\mu^{\pm}(x)= \mu^{\pm}(x,W^{\pm}) := \frac{\exp(-W^{\pm}(x-1))+\exp(-W^{\pm}(x))}{2\sum_{x \in \Z}\exp(-W^{\pm}(x))}.
\end{eqnarray}
Let $Q^{+}$ (resp. $Q^{-}$) the distributions of  $\mu^{+}$ (resp. $\mu^{-}$), define $Q=1/2(Q^++Q^-)$ the distribution of a random measure called $\mu$. Let $l^1=\{l:\Z \rightarrow \R,\ ||l||:=\sum_{x \in \Z} |l(x)|<+\infty\}$, a main ingredient in the proof of the theorems is  the following result (\cite{GanPerShi} Corollary 1.1), for all functions, $f:l^1\rightarrow \R$ which is shift-invariant,  
\begin{eqnarray}
\lim_{n \rightarrow +Ê\infty}\E\left[f\left(\frac{\lo(x,n)}{n}, x \in \Z\right)\right]= \E[f(\{\mu(x), x \in \Z\})].  \label{1}
\end{eqnarray}

\subsection{Proof of Theorem \ref{th1}}
Define 
\begin{eqnarray*}
R_n(r)= \sup_{x\in \Z} \sum_{k \in I_{r}(x)}\frac{\lo(k,n)}{n}, 
\end{eqnarray*}
where $I_{r}(x)$ is defined above \ref{eq1}, from \ref{1} we get
\begin{eqnarray}
R_n(r)\rightarrow \sup_{x\in \Z} \sum_{k \in I_{r}(x)}\mu(x), \textrm{ in law}. \label{2.2}
\end{eqnarray}
 Moreover it was proven in \cite{Pierre2} (see also \cite{GanShi}) that
\begin{eqnarray}
\limsup_n R_n(r)= \textrm{const} \in  [0,+\infty), \ \p. a. s., \label{11}
\end{eqnarray}
then \ref{2.2} and \ref{11} yield the, easy to get fact
\begin{eqnarray}
\limsup_n R_n(r) \geq g(r) := \sup\left \{z,\ z \in supp \left( \sup_{x\in \Z} \sum_{k \in I_{r}(x)}\mu(k) \right) \right\}, \ \p. a. s.  \label{2.6b}
\end{eqnarray}
The next step is to determine which environment maximizes $\sup_{x\in \Z} \sum_{k \in I_{r}(x)}\mu(k)$. It is not surprising that this environment is almost the same of the media that maximizes $\sup\left \{z,\ z \in supp \left( \sup_{x\in \Z} \mu(x) \right) \right\}$ in \cite{GanPerShi}. We will restrict our analysis to $\mu^+$ (the restriction of $\mu$ to $Q^+$), and shortly discuss what happens under $Q^-$. \\
The definition of the maximizing environment (\textit{the infinite valley}) is given by : for all $x \in \Z^*_+$ assume $\alpha_x^+:=\bar{\alpha}$, for all $x \in \Z^*_-$, $\alpha_x^+:=\bar{A}$ and $\alpha_0^+=\bar{\alpha}\un_{\{1-\bar{A} \geq \bar{\alpha}\}}+\bar{A}\un_{\{1-\bar{A} < \bar{\alpha}\}}\equiv \bar{\alpha}\un_{\{\tilde{A} \geq \tilde{\alpha}\}}+\bar{A}\un_{\{\tilde{A} < \tilde{\alpha}\}} $. Notice that in \cite{GanPerShi}, the value of $\alpha_0$ is chosen in such a way that the random walk starting from $0$ is pushed on the site (+1 or -1) from where it as a better chance to get back to 0. Here things are a bit different, and the value of $\alpha_0$ is taken in such a way that the walk starting from 0 is pushed to the side where it will encounter the center of the interval of concentration.

With this choice we have the following expression for the exponential of the maximizing potential $\bar{W}^+$:
\begin{eqnarray}
\exp(-\bar{W}^+(x))=\left\{ \begin{array}{ll} \tilde{\alpha}^x & \textrm{if}\ x>0, \\
1 & \textrm{if }  x=0, \\
   \tilde{A}^{-x-1} \left(\frac{1}{\tilde{\alpha}}\un_{\{\tilde{\alpha}<\tilde{A}\}}+\tilde{A}\un_{\{\tilde{\alpha}\geq \tilde{A}\}}\right) & \textrm{if }  x <0.     \end{array} \right. \label{2.6}
 \end{eqnarray}
From \ref{2.6} we easily get the following expression for $g(r)$:
\begin{eqnarray}
& & g(r)=\frac{(1-\tilde{\alpha})(1-\tilde{A})}{2(1-\tilde{\alpha}\tilde{A})}\sup_{x\in \Z} \left(\sum_{k \in I_r(x)} F(k)+F(k-1)\right), \label{eq2.7} \textrm{ where we recall } \\
& & F(l)=\left\{ \begin{array}{ll} \tilde{\alpha}^l(\tilde{\alpha}\un_{\{\tilde{\alpha}<\tilde{A}\}}+\un_{\{\tilde{\alpha} \geq\tilde{A}\}}) & \textrm{if}\ l>0, \\
\tilde{\alpha}\un_{\{\tilde{\alpha}<\tilde{A}\}}+\un_{\{\tilde{\alpha} \geq\tilde{A}\}} & \textrm{if }  l=0, \\
   \tilde{A}^{-l-1} (\un_{\{\tilde{\alpha}<\tilde{A}\}}+\tilde{A}\un_{\{\tilde{\alpha} \geq\tilde{A}\}}) & \textrm{if }  l <0.     \end{array} \right.
\end{eqnarray}
of course the value of the supremum in  \ref{eq2.7} depends on $\tilde{A}$ and $\tilde{\alpha}$. Let us denote $a^+ \in \Z$ the coordinate 
of this supremum. Notice that $a^+$ is deterministic but not necessarily equal to $0$, moreover  if $a^+<0$ then $0 \in I_r(a^+)$ whereas if $a^+>0$ either $0 \in I_r(a^+)$ or $1= \min\{x,\ x \in  I_r(a^+)\}$. This because the maximum of $F$ is 1.

To turn the inequality \ref{2.6b} into an equality  we have to show that 
$$\limsup_n R_n(r) \leq g(r),\  Q^+\otimes \pa.a.s.$$
To get this last inequality we use the same method as in \cite{GanPerShi}. The proof is based on two facts, the first one says that the distribution of $(\lo(I_r(x),n),n\in \N)$ under $\pa$ is stochasticaly dominated by the distribution of $(\lo(I_r(a^+),n),n\in \N)$ under $\p^{\bar{W}^+}$ (Lemma 3.2 in \cite{GanPerShi}). The second fact says that, under $\p^{\bar{W}^+}$ $(X_n,n\in \N)$ is a positive recurrent Markov chain with invariant probability measure given by $\mu^+$, therefore, as $I_r(a^+)$ is finite, for all $\epsilon>0$, $\p^{\bar{W}^+}\left[\sum_{k \in I_r(a^+)}\frac{\lo(k,n)}{n} \geq \sum_{k \in I_r(a^+)}\bar{\mu}^+(x)+\epsilon \right]\leq \exp(-C(\epsilon)n)$. $C(\epsilon)$ is strictly positive constant depending only on $\epsilon$. We finally get that 
 \begin{eqnarray}
\limsup_n R_n(r) =g(r), \ Q^+\otimes \pa. a. s. \label{eq5}
\end{eqnarray}
 If we do the same computations with respect to $Q^-$, the above result remains the same, indeed we can take the same definition for $\bar{W}^-$ we took for $\bar{W}^+$ because both are  equal to $0$ only in $0$.
In order to get Theorem \ref{th1} from \ref{eq5} it suffices to notice that $Y_{n,g(r)}$ is in some sense the dual of $ \sup_{x \in \Z}\sum_{k \in I_{r}(x)}\frac{\lo(k,n)}{n}$, indeed, let $0<\beta<1$, $\gamma>0$ we easily get the following assertions,
\begin{eqnarray} 
 \limsup_n R_n(f(\beta)) \leq \beta &\Rightarrow & \liminf_n Y_{n, \beta} \geq f(\beta)-\gamma, \label{2.10} \\
 \limsup_n R_n(f(\beta)) \geq \beta &\Rightarrow & \liminf_n Y_{n, \beta-\gamma} \leq f(\beta), \label{2.11}
\end{eqnarray} 
where $f(\beta)$ is the smallest integer that satisfies the following inequality: 
 \begin{eqnarray}
 g(f(\beta))\geq \beta.
 \end{eqnarray}
 By collecting the last four equations and by letting $\gamma$ goes to 0 ($\liminf_n Y_{n, \beta-\gamma}$ increases when $\gamma \searrow 0$), we finally get that $\liminf Y_{n}(\beta)=f(\beta)\ \p.a.s.$ 
 
\subsection{Proof of Theorem \ref{th2}}
\begin{figure}[h]
\begin{center}
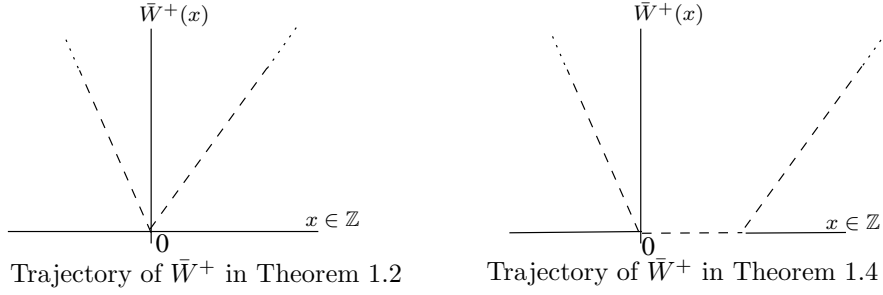 
\caption{Infinite valleys, $(1-\bar{\alpha}<\bar{A})$} \label{fig7}
\end{center}
\end{figure}

First we notice that for all fixed $\delta$, $Z_{n,\delta}$ is bounded because $\lo(k,n)\leq n$, in fact it is also clear that $Z_{n,\delta}\leq 1/\delta$, so we can apply the 3 steps of the proof of Lemma 3 of \cite{GanShi} (see also Section 2.1 in \cite{Pierre2})  for $Z_{n,\delta}$ instead of $\lo^*(n)/n$ to show that $\limsup_n Z_{n,\delta}=\textrm{const}\in [0,+\infty),\ \p.a.s.$, with this result and \ref{1} we easily get that 
\begin{eqnarray} \limsup_n \sum_{x \in \Z}\un_{ \{\lo(k,n)/n \geq \delta \}} \geq m(\delta)=m(\delta,\mu):= \sup\left \{z,\ z \in supp \left(\sum_{x \in \Z}\un_{\{{\mu(x)\geq \delta} \}}\right) \right\}, \ \p. a. s.
\end{eqnarray}
As before we need to determine a random environment that maximizes $\sum_{x \in \Z}\un_{\{{\mu(x)\geq \delta} \}}$, this time it is not exactly the same as the preceding one. Indeed we build an  infinitely deep  valley such that the asymptotic of the normalized local time reach the level $\delta$ and also flat at the bottom of this valley (see Figure \ref{fig7}) such that there is a large number of points where the normalized local time reach this level $\delta$. We easily get that such an environment leads us to
\begin{eqnarray*}
\exp(-\bar{W}^+(x))=\left\{ \begin{array}{ll} \tilde{\alpha}^{x} & \textrm{if }\ x > g(\delta), \\
1 & \textrm{if } 0 \leq x \leq g(\delta),  \\
   \tilde{A}^{-x} & \textrm{if }  x <0,    \end{array} \right.  \ \textrm{under }Q^+
 \end{eqnarray*}
and 
\begin{eqnarray*}
\exp(-\bar{W}^-(x))=\left\{ \begin{array}{ll} \tilde{\alpha}^x & \textrm{if }\ x > 0, \\
1 & \textrm{if } -g(\delta) \leq x \leq 0,  \\
   \tilde{A}^{-x-1}\frac{1}{\tilde{\alpha}} & \textrm{if }  x <-g(\delta).    \end{array} \right.  \textrm{under }Q^-
 \end{eqnarray*}
where $g(\delta)$ is a positive integer that will be determined later. Under $Q^+$, we get that
\begin{eqnarray}
\bar{\mu}^+(x):=\mu^+(x,\bar{W}^+)= \frac{\exp(-\bar{W}^+(x-1))+\exp(-\bar{W}^+(x))}{2(\tilde{\alpha}^{g(\delta)+1}/(1-\tilde{\alpha})+g(\delta)+1+\tilde{A}/(1-\tilde{A}))}
 \end{eqnarray}
and a similar expression for $\bar{\mu}^-:=\mu^-(x,\bar{W}^-)$. Notice that for all $x\in G(\delta):= \{1,\cdots,g(\delta)\}$, $\bar{\mu}^+(x)$ reaches its maximum which is equal to $1/(\bar{\alpha}^{g(\delta)+1}/(1-\bar{\alpha})+g(\delta)+1+\bar{A}/(1-\bar{A}))$.
Let us decompose $\bar{m}(\delta):=m(\delta,\bar{\mu}^+)$ in the following way
\begin{eqnarray*}
& & \bar{m}(\delta):=\bar{m}_S(\delta)+\bar{m}_L(\delta), \textrm{ with}\\
& & \bar{m}_S(\delta):=\sum_{x \in \Z \smallsetminus G(\delta) }\un_{\{{\bar{\mu}^+(x)\geq \delta} \}}, \textrm{ and}\\
& & \bar{m}_L(\delta):=\sum_{x \in G(\delta) }\un_{\{{\bar{\mu}^+(x)\geq \delta} \}}.
\end{eqnarray*}
We easily get that
\begin{eqnarray*}
& & \bar{m}_L(\delta)=0 \Rightarrow \bar{m}_S(\delta)=0, \\
& & \bar{m}_L(\delta) >0 \Rightarrow \bar{m}_L(\delta)=g(\delta), 
\end{eqnarray*}
therefore 
\begin{eqnarray}
\bar{m}(\delta)\un_{\{\bar{m}(\delta)>0\}}=g(\delta)+\bar{m}_S(\delta) \geq g(\delta). 
\end{eqnarray}
Moreover to get $\bar{m}(\delta)$ as large as possible,  $g(\delta)$ must be the largest integer that satisfies the following
$$1/(\tilde{\alpha}^{g(\delta)+1}/(1-\tilde{\alpha})+g(\delta)+1+\tilde{A}/(1-\tilde{A}))\geq \delta, $$ and therefore the largest integer such that
$$g(\delta) \leq \frac{1}{\delta}-1-\frac{\tilde{\alpha}^{g(\delta)+1}}{1-\tilde{\alpha}}-\frac{\tilde{A}}{1-\tilde{A}}.$$
\noindent The same computations lead to a similar expression under $Q^-$ ($g(\delta) \leq \frac{1}{\delta}-1-\frac{1}{\tilde{\alpha}}\frac{\tilde{A}^{g(\delta)}}{1-\tilde{A}}-\frac{\tilde{\alpha}}{1-\tilde{\alpha}}$). Assembling what we did above yields to 
\begin{eqnarray*}  \frac{1}{\delta}-1-\max\left(\frac{1}{\tilde{\alpha}}\frac{\tilde{A}^{g(\delta)}}{1-\tilde{A}}+\frac{\tilde{\alpha}}{1-\tilde{\alpha}},\frac{\tilde{\alpha}^{g(\delta)+1}}{1-\tilde{\alpha}}+\frac{\tilde{A}}{1-\tilde{A}}\right) \leq \limsup_n Z_{n,\delta} \leq \frac{1}{\delta},\ \p. a. s.
\end{eqnarray*}
As $\tilde{A}<1$ and $\tilde{\alpha}<1$, we get the theorem.

\bibliographystyle{alpha}

 \bibliography{thbiblio}

\vspace{1cm} \noindent
\begin{tabular}{l}
Laboratoire MAPMO - C.N.R.S. UMR 6628 - F\'ed\'eration Denis-Poisson  \\
Universit\'e d'Orl\'eans, UFR Sciences \\
B\^atiment de math\'ematiques - Route de Chartres \\
B.P. 6759 - 45067 Orl\'eans cedex 2 \\
FRANCE 
\end{tabular}

\end{document}